\documentclass{article}
\usepackage[russian]{babel}
\usepackage[cp1251]{inputenc}
\usepackage[T2B]{fontenc}
\usepackage[pdftex, colorlinks, urlcolor=blue, pdfborder={0 0 0 [2 0]}]{hyperref}
\usepackage{amsfonts,amsmath,amssymb}
\usepackage{amsthm}
\usepackage{array}
\usepackage{indentfirst}

\newtheorem{theorem}{Теорема}

\newtheorem{proper}{Утверждение}

\newtheorem{m-lemma}{Микролемма}
\newenvironment{proof-m}
{\vspace{1pt}\par{\sl%
Д\,о\,к\,а\,з\,а\,т\,е\,л\,ь\,с\,т\,в\,о\ \ м\,и\,к\,р\,о\,л\,е\,м\,м\,ы.\,\ }}%
{\noindent\vspace{1pt}}

\newenvironment{proofs}
{\vspace{1pt}\par{\sl%
Д\,о\,к\,а\,з\,а\,т\,е\,л\,ь\,с\,т\,в\,о\ \ т\,е\,о\,р\,е\,м\,ы.\,\ }}%
{\noindent\vspace{1pt}}

\newenvironment{proof-p}
{\vspace{1pt}\par{\sl%
Д\,о\,к\,а\,з\,а\,т\,е\,л\,ь\,с\,т\,в\,о\ \ у\,т\,в\,е\,р\,ж\,д\,е\,н\,и\,я.\,\ }}%
{\noindent\vspace{1pt}}

\makeatletter
\renewcommand{\@biblabel}[1]{#1.} 
\makeatother

\begin{document}

\title{"<Прямой"> метод доказательства обобщенной формулы Ито -- Вентцеля для обобщенного стохастического дифференциального уравнения}

\author{Карачанская  Е.В. \\
{Тихоокеанский госуниверситет, Хабаровск}}



\date{}

\maketitle
\begin{abstract}
В данной статье представлено полностью одно из доказательств обобщенной формулы Ито -- Вентцеля, идеи которого изложены в  \cite{12_KchDubPrep} ({\it arXiv}:1309.3038v1). Данное доказательство представляет подход к построению обобщенной формулы Ито -- Вентцеля, основанный на непосредственном привлечении обобщенной формулы Ито и теории стохастической аппроксимации.
\end{abstract}

Ключевые слова: формула Ито -- Вентцеля, обобщенное уравнение Ито, нецентрированная пуассоновская мера, аппроксимация $\delta-$функции.

\section*{Введение}
Важным элементом теории стохастических дифференциальных уравнений является построение дифференциалов от функций, зависящих от случайного процесса, который подчинен стохастическим уравнениям. Такими являются формула Ито \cite{Ito_59,Kunita-Watanabe-rus_71} -- дифференциал  неслучайной функции от случайного процесса; формула Ито -- Вентцеля \cite{Wentzel_65}, позволяющая построить дифференциал для функции, которая сама является решением стохастического уравнения. Множество работ посвящено этому вопросу \cite{Rozovsky_73,D_89,Ocone_89} и др. Естественно было рассмотреть не только случай наличия винеровских процессов в СДУ Ито, но и когда присутствует пуассоновский процесс. В 2002 г. на основе представления о ядрах интегральных инвариантов в работе В.~А.~Дубко  \cite{D_02} было построено обобщение стохастических дифференциалов от случайных функций, подчиненных обобщенному уравнению Ито (ОСДУ), которое было названо им {\it обобщенным уравнением Ито -- Вентцеля}. То есть в отличие от других работ, результат связан не с расширением формулы Ито -- Вентцеля на более широкий класс функций, а получена новая формула, связанная с обобщенным уравнением Ито. В 2007 г. наиболее удачная с нашей точки зрения попытка получения обобщенной формулы Ито -- Вентцеля на основе классической теории СДУ Ито (который будем называть в дальнейшем прямым методом) была предпринята  B.~Oksendal и  T.~Zhang \cite{Oks_07}, при отсутствии винеровской составляющей для одномерного процесса. В 2011 г. в работе  \cite{11_KchOboz},  снова опираясь на уравнения для ядер интегральных инвариантов, была предложена обобщенная формула Ито -- Вентцеля для нецентрированной пуассоновской меры (НЦПМ) -- в отличие от работы \cite{D_02}. Отметим, что хотя теоремы существования и единственности опираются на центрированные пуассоновские меры (ЦПМ), для целей программного управления более адекватным является применение уравнений с нецентрированными мерами \cite{11_KchUpr}. Заметим также, что в этом случае требование о характере пуассоновского распределения носит только общее ограничение, не требующее знания явного вида. В то же время в задачах программного управления, если опираться на центрированную меру, требуется явный, в не общий,  вид распределений интенсивности пуассоновских скачков \cite{D_02}.
По-видимому, с нецентрированной мерой результат получен  впервые в работе \cite{11_KchOboz}.

Цель данной работы -- представить метод получения обобщенного уравнения Ито -- Вентцеля для случая НЦПМ,  опираясь на классические положения теории случайных процессов и теории стохастической аппроксимации.

 В работе предлагается в качестве индикатора использовать $\delta$-функцию, рассматривая ее как предел последовательности экспонент. Поскольку экспонента -- бесконечно дифференцируемая функция (нам достаточно второй производной), и соответствующая последовательность быстро сходится, и такое допредельное представление позволяет применять обобщенную формулу Ито без дополнительных замечаний.
Пример  формального использования $\delta $-функции при построении  стохастических аналогов задач  классической механики на основе уравнений Ито был приведен в работе~\cite{D_84}.

Отметим, что представленная в \cite{11_KchOboz} модификация обобщенной формулы Ито -- Вентцеля, полученная на основе ядер интегральных инвариантов, требует более жестких условий для коэффициентов всех рассматриваемых уравнений (существование вторых производных) \cite{12_KchDubPrep}.

\section{Предварительные замечания: эволюция формул дифференцирования сложной функции}

Пусть ${\bf x}(t)\in \mathbb{R}^{n} $ решение уравнения
\begin{equation*}
\frac{d{\bf x}(t)}{dt} =A(t),\ \ \ \ {\bf x}(t)|_{t=0} ={\bf x}(0).
\end{equation*}
Пусть функция $F(t;{\bf x})$ непрерывна вместе со своими частными производными $F'_{t} (t;{\bf x})$, $F'_{x_{i} } (t;{\bf x})$, $i=\overline{1,n}$. Тогда, опираясь на  правило дифференцирования сложной функции, получаем:
\begin{equation*}
\frac{dF(t;{\bf x}(t))}{dt} =F'_{t} (t;{\bf x}(t))+\sum\limits_{j=1}^{n}a_{j} (t)F'_{x_{j} } (t;{\bf x}(t)),\ \ \ \ F(t;{\bf x})|_{t=0} =F(0;{\bf x}).
\end{equation*}
Если же $F(t;{\bf x})$ есть решение уравнения
\begin{equation*}
\frac{\partial F(t;{\bf x})}{\partial t} =Q(t;{\bf x}), \ \ \ \ \ F(t;{\bf x})|_{t=0} =F(0;{\bf x}) ,
\end{equation*}
то переходим к выражению:
\begin{equation}\label{odu-prav}
\frac{dF(t;{\bf x}(t))}{dt} =Q(t;{\bf x}(t))+\sum\limits_{j=1}^{n}a_{j} (t)F'_{x_{j} } (t;{\bf x}(t)).
\end{equation}
Рассмотрим случай, когда ${\bf x}(t)$ определяется уравнением
\begin{equation} \label{GrindEQ__1_1_2_}
 d{\bf x}(t)=A(t;{\bf x})dt+\sum\limits_{k=1}^{m}B_{k} (t;{\bf x})dw_{k} (t).
\end{equation}
Если $F(t;{\bf x})$ -- случайная функция, то ее дифференциал определяется формулой, называемой \textit{формулой Ито} \cite{Ito_59}. Далее, пусть  $F(t;{\bf x})$ подчинено уравнению
\begin{equation} \label{GrindEQ__1_4_1_} dF(t;{\bf x})=Q(t;{\bf x})dt+\sum\limits_{k=1}^{m}D_{k} (t;{\bf x})dw_{k} (t),
\end{equation}
где ${\bf w}(t)$ -- тот же  $m$-мерный винеровский процесс с  независимыми  компонентами, как и в уравнении \eqref{GrindEQ__1_1_2_}.
Относительно случайных функций $A(t;{\bf x})$, $B_{k}(t;{\bf x})$, $Q(t;{\bf x})$ и $D_{k}(t;{\bf x})$ предполагаем, что они  непрерывны и ограниченные вместе со своими первыми  и вторыми частными производными по компонентам ${\bf x}$, неупреждающие по $t$ относительно приращений винеровского векторного процесса ${\bf w}(t)$ и согласованные с потоком $\sigma$-алгебр, индуцированным ${\bf w}(t)$ на всем интервале $[0,T]$. Этих ограничений достаточно для построения  дифференциала от процесса $F(t;{\bf x}(t))$:
\begin{equation} \label{GrindEQ__1_4_2_}
\begin{array}{c}
\displaystyle d_{t} F(t;{\bf x}(t))=Q(t;{\bf x}(t))dt+\sum\limits_{k=1}^{m}D_{k} (t;{\bf x}(t))dw_{k} +\\
+\displaystyle \Bigl[\sum\limits_{i=1}^{n}a_{i} (t)\frac{\partial F(t;{\bf x})}{\partial x_{i} }\bigr|_{{\bf x}={\bf x}(t)} +\frac{1}{2}\sum\limits_{i=1}^{n} \sum\limits_{j=1}^{n}\sum\limits_{k=1}^{m}b_{i,k} (t)b_{j,k} (t)\frac{\partial ^{{\kern 1pt} 2} F(t;{\bf x})}{\partial x_{i} \partial x_{j} }\bigr|_{{\bf x}={\bf x}(t)}\Bigr. + \\
+\Bigl.\displaystyle
\sum\limits_{i=1}^{n}\sum\limits_{k=1}^{m}b_{i,k} (t)\frac{\partial D_{k} (t;{\bf x})}{\partial x_{i} }\bigr|_{{\bf x}={\bf x}(t)} \Bigr]dt+
\sum\limits_{i=1}^{n}b_{i,k} (t)\displaystyle\frac{\partial F(t;{\bf x})}{\partial x_{i} }\bigr|_{{\bf x}={\bf x}(t)} dw_{k}, \end{array}
 \end{equation}
\begin{equation*}
F(t;{\bf x}(t))\left|_{t=0} =F(0;y),\right. \; F(0;y)\in C_{0}^{2} .
\end{equation*}
Соотношение  \eqref{GrindEQ__1_4_2_} носит название \textit{формула Ито -- Вентцеля} \cite{Wentzel_65,Rozovsky_73}.

Перейдем к построению подобной формулы  для случая,  когда случайный процесс ${\bf x}(t)$ является решением уравнения
\begin{equation}\label{1}
d{\bf x}(t)=A(t)dt+B(t)d{\bf w}(t)+\int g(t;\gamma) {\rm \; }\nu (dt;d\gamma ),
\end{equation}
где ${\bf x}(t) \in \mathbb{R}^{n} $; ${\bf w}(t)$ -- $m$-мерный винеровский процесс; $\nu (dt;d\gamma )$ -- нецентрированная мера Пуассона.

Для случайной функции $F(t;{\bf x}(t))$, где ${\bf x}(t)$ --
решение уравнения (\ref{GrindEQ__2_1_1_}), дифференциал представлен \textit{обобщенной формулой Ито} \cite{Kunita-Watanabe-rus_71}
(используем в виде \cite[с.\, 271-272]{GS_68}):

\begin{equation}\label{Ayd2}
\begin{array}{c}
 \displaystyle d_{t}F(t;{\bf x}(t))=\Bigl[ \frac{\partial F(t;{\bf x}) }{\partial t}\bigl|_{{\bf x}={\bf x}(t)}+
 \sum\limits_{i=1}^{n}a_{i}(t;{\bf x}(t))\frac{\partial F(t;{\bf x}) }{\partial x_{i}}\bigl|_{{\bf x}={\bf x}(t)}+\Bigr.\\
\Bigl.+ \displaystyle  \frac{1}{2} \sum\limits_{i=1}^{n}
 \sum\limits_{j=1}^{n}
 \sum\limits_{k=1}^{m}b_{i\,k}(t;{\bf x}(t))b_{j\,k}(t;{\bf x}(t))
 \frac{\partial^{\,2} F(t;{\bf x}) }{\partial x_{i}\partial x_{j}}\bigl|_{{\bf x}={\bf x}(t)}\Bigr]dt +\\
 +\displaystyle \sum\limits_{i=1}^{n}
 \sum\limits_{k=1}^{m}b_{i\,k}(t;{\bf x}(t))\frac{\partial F(t;{\bf x}) }{\partial x_{i}}\bigl|_{{\bf x}={\bf x}(t)} dw_{k}(t)+\\
 +\displaystyle \int
 \Bigl[ F\left(t;{\bf x}(t)+g(t;{\bf x}(t);\gamma)\right)- F(t;{\bf x}(t))
 \Bigr]\nu(dt;d\gamma).
 \end{array}
\end{equation}

Получим правило дифференцирования для $F(t;{\bf x}(t))$, когда $F(t;{\bf x})$ -- решение уравнения
\begin{equation*}
dF(t;{\bf x})=Q(t;{\bf x})dt+\sum\limits_{k=1}^{m}D_{k} (t;{\bf x})dw_{k} (t)+\int  G(t;x;\gamma ) \nu (dt;d\gamma ),
\end{equation*}
а ${\bf x}(t)$ -- решение уравнения \eqref{1}.

Следуя традиции, полученную в дальнейшем формулу будем называть \textit{обобщенной формулой Ито -- Вентцеля}.

Первоначально эта формула была получена именно на основе уравнений для ядер интегральных инвариантов,  в связи с поиском сохраняющихся функционалов для открытых систем \cite{D_02,11_KchOboz}. Попытки получить эту формулу прямым методом предпринималась в работе \cite{Oks_07}. Именно последняя работа и стимулировала нас к рассмотрению общего случая: присутствие винеровских и пуассоновских возмущений и в $n$-мерном случае для процесса ${\bf x}(t)$ \cite{12_KchDubPrep}.

Особенность представленного подхода состоит в том, что в отличие  от известных работ других авторов, мы с самого начала будем работать c обобщенными уравнениями для нецентрированных мер. Это требует несколько более жестких ограничений на класс весовых функций при пуассоновских возмущениях, но делает доказательство более компактным и прозрачным.

Полученная ранее модификация обобщенной формулы Ито -- Вентцеля исходила из физических соображений (ядра интегральных инвариантов). Другой ее вариант докажем "<прямым"> методом, основываясь на классических представлениях, при этом будем рассматривать более общий случай, связанный с видом коэффициентов, хотя она получается при менее жестких ограничениях, чем через ядра. Связано это с тем, что ядра существуют при определенных ограничениях на коэффициенты.

\section{Основной результат}
Будем использовать обозначения, аналогичные \cite{GS_68}.
$H_{n}$ -- означает пространство случайных векторов (отображений) $\mathcal{F}_t$-измеримых при произвольном $t\in [0,T]$ и с вероятностью 1  и таких, что (\cite[c.\,255]{GS_68}):
\begin{equation*}
\int\limits_{0}^{T}|\alpha (t)|^{n} dt<\infty.
\end{equation*}
Если \cite[c.\,257]{GS_68} функция $\varphi (t;y)=\varphi (t;y;\omega )$, при произвольных $t\in [0,T]$ и $y\in \mathbb{R}$, измеримая как функция трех аргументов $t;y;\omega $, и
\begin{equation*}
\int\limits_{0}^{T}\int\limits_{{\mathbb R}}  |\varphi (t;y)|^{2} \Pi (dy)dt<\infty  \end{equation*}
то используется обозначение: $\varphi (t;y)\in H_{2} (\Pi)$.

Рассмотрим случайный процесс ${\bf x}(t)$ со значениями в $\mathbb{R}^{n}$, определяемый следующим уравнением:
\begin{equation}\label{GrindEQ__2_1_1_}
d{\bf x}(t)=A(t)dt+B(t)d{\bf w}(t)+\int g(t;\gamma )\nu (dt;d\gamma),
\end{equation}
где $A(t)=\{ a_{1} (t),...,a_{n} (t)\}^{*}$, $g(t;\gamma )=\{ g_{1} (t;\gamma ),...,g_{n} (t;\gamma )\}^{*} \in \mathbb{R}^{n} $,  $\gamma\in \mathbb{R}^{n'} $,  ${\bf w}(t)$ -- $m$-мерный винеровский процесс;
\begin{equation} \label{GrindEQ__2_1_2_}
\begin{array}{c}
  |g(t;\gamma )|\in H_{1,2} (\Pi ),\ \ \ \ \sqrt{|a_{j} (t)|},\, |b_{j,k} (t)|\in H_{2} , \\
  B(t)=\{ b_{j,k} (t); \ \ \ j=\overline{{1,n}};\ \ \ k=\overline{{1,m}},\}
\end{array}
\end{equation}
$\nu (\Delta t;\Delta \gamma )$ - стандартная пуассоновская мера на $[0,T] \times \mathbb{R}^{n'}$, моделирующая независимые случайные величины на непересекающихся интервалах и множествах.

\begin{theorem} \label{t1}
Пусть $F(t;{\bf x})$, $(t;{\bf x})\in[0;T]\times\mathbb{R}^{n} $ -- скалярная функция, обобщенный стохастический дифференциал которой имеет вид:
\begin{equation} \label{GrindEQ__2_5_1_}
d_{t} F(t;{\bf x})=Q(t;{\bf x})dt+\sum\limits_{k=1}^{m}D_{k} (t;{\bf x})dw_{k} (t)+\int G(t;{\bf x};\gamma ) \nu (dt;d\gamma )
\end{equation}
и для коэффициентов \eqref{GrindEQ__2_5_1_} выполнены условия:
\begin{description}
  \item[$\mathfrak{L.}a)$:] $Q(t;{\bf x})$, $D_{k} (t;{\bf x})$, $G(t;{\bf x};\gamma )\in\mathbb{R}$ в общем случайные функции, измеримые относительно того же неубывающего  потока $\sigma$-алгебр $\bigl\{\mathcal{F}_{t}\bigr\}_{0}^{T}$, что и процессы $w(t)$, $\nu (t;\mathcal{A})$ для любого множества  $\mathcal{A}\in \mathfrak{B}$  из фиксированной борелевской $\sigma $-алгебры  {\rm(\cite[с. 266]{GS_68})};
  \item[$\mathfrak{L.}b)$:]c вероятностью единица  функции $Q(t;{\bf x}),D_{k} (t;{\bf x}),G(t;{\bf x};u)$ непрерывны и ограничены по совокупности переменных вместе со своими вторыми частными производными по компонентам вектора ${\bf x}$.
\end{description}
Тогда, если случайный процесс ${\bf x}(t)$  подчинен системе  \eqref{GrindEQ__2_1_1_} и  выполняются ограничения \eqref{GrindEQ__2_1_2_}, то существует стохастический дифференциал:
\begin{equation} \label{GrindEQ__2_5_2_}
\begin{array}{c}
 \displaystyle d_{t} F(t;{\bf x}(t))=Q(t;{\bf x}(t))dt+\sum\limits_{k=1}^{m}D_{k} (t;{\bf x}(t))dw_{k} +\\
+\displaystyle \Bigl[\sum\limits_{i=1}^{n}a_{i} (t)\frac{\partial F(t;{\bf x})}{\partial x_{i} }\bigr|_{{\bf x}={\bf x}(t)} +\frac{1}{2} \sum\limits_{i=1}^{n}\sum\limits_{j=1}^{n}\sum\limits_{k=1}^{m}b_{i,k} (t)b_{j,k} (t)\frac{\partial ^{{\kern 1pt} 2} F(t;{\bf x})}{\partial x_{i} \partial x_{j} }\bigr|_{{\bf x}={\bf x}(t)}\Bigr. + \\
+\Bigl.\displaystyle
\sum\limits_{i=1}^{n}b_{i,k} (t)\frac{\partial D_{k} (t;{\bf x})}{\partial x_{i} }\bigr|_{{\bf x}={\bf x}(t)} \Bigr]dt+
\sum\limits_{i=1}^{n}\sum\limits_{k=1}^{m}b_{i,k} (t)\displaystyle\frac{\partial F(t;{\bf x})}{\partial x_{i} }\bigr|_{{\bf x}={\bf x}(t)} dw_{k} + \\
+\displaystyle\int \Bigl[\bigl(F(t;{\bf x}(t)+g(t;\gamma )\bigr)-F(t;{\bf x}(t))\Bigr]\nu (dt;d\gamma ) +\\
+ \displaystyle\int G\bigl(t;{\bf x}(t)+g(t;\gamma );\gamma\bigr)\nu (dt;d\gamma ). \end{array}
\end{equation}
\end{theorem}

Выражение \eqref{GrindEQ__2_5_2_} -- \textit{обобщеннная формула Ито -- Вентцеля} на случай присутствия как винеровских, так и пуассоновских возмущений.

\section{"<Прямой"> метод доказательства}

Доказательство базируется на определении $\delta $ - функции как предела некоторой специальной функции.

\begin{proper}\label{lem-eps}
Пусть функция $f(t;x)$ ограничена с вероятностью единица для любых $t\in [0,T]$  и  удовлетворяет условию Липшица по компоненте~$x$. Тогда имеет место представление
\begin{equation} \label{GrindEQ__2_5_3_}
\begin{array}{l}
f(t;x)=\displaystyle\int\limits_{-\infty }^{\infty }f(t;y) \delta (y-x)dy=
\displaystyle\lim\limits_{\varepsilon \downarrow 0} \int\limits_{-\infty }^{\infty }f(t;y) \delta_{\varepsilon } (y-x) dy,
\end{array}
\end{equation}
где $\delta_{\varepsilon } (y-x)
 =\displaystyle\frac{1}{\varepsilon \sqrt{2\pi } } \exp \left\{ -\frac{(y-x)^{2} }{2\varepsilon ^{2} } \right\} $.
\end{proper}
\begin{proof-p}
Положим, что для функции $f(t;x)$ выполняется условие Липшица:
\begin{equation*}|f(t;y_{1} )-f(t;y_{2} )|\le L(t)|y_{1} -y_{2} |^{\varsigma } ,0<\varsigma \le 1.
\end{equation*}
Поскольку $t$ в доказательствах рассматривается как фиксированное, то для упрощения в записи формул, параметр $t$ будем опускать. Т.е., вместо $f(t;x)$ будем использовать $f(x)$, а  вместо $L(t)$ -- обозначение  $L$.

Выполним замену переменных
$(y-x)\varepsilon ^{-1} =z,{\rm \; \; \; }y=\varepsilon z+x .
$
Добавив и вычитая одно и то же выражение, и, используя свойства интеграла вероятности, получаем:
\begin{equation}\label{dob1}
\begin{array}{c}
  \displaystyle\frac{1}{\sqrt{2\pi } } \int\limits_{-\infty }^{\infty }f(\varepsilon z+x) \exp \left\{ -\frac{z^{2} }{2} \right\} dz=f(x)\frac{1}{\sqrt{2\pi } } \int\limits_{-\infty }^{\infty } \exp \left\{ -\frac{z^{2} }{2} \right\} dz+
  \\
  +\displaystyle\frac{1}{\sqrt{2\pi } } \int\limits_{-\infty }^{\infty }\left(f(\varepsilon z+x) -f(x)\right)\exp \left\{ -\frac{z^{2} }{2} \right\} dz=\\
  = f(x) +
\displaystyle\frac{1}{\sqrt{2\pi } } \int\limits_{-\infty }^{\infty }\left(f(\varepsilon z+x) -f(x)\right)\exp \left\{ -\frac{z^{2} }{2} \right\} dz
\end{array}
\end{equation}
Оценим последний интеграл:
\begin{equation}\label{ocenka1}
\begin{array}{c}
  \displaystyle\left|\frac{1}{\sqrt{2\pi } } \int\limits_{-\infty }^{\infty }\left(f(\varepsilon z+x) -f(x)\right)\exp \left\{ -\frac{z^{2} }{2} \right\} dz\right|\le \\
 \le \displaystyle\frac{L}{\sqrt{2\pi } } \int\limits_{-\infty }^{\infty }|(\varepsilon z+x) -x|^{\varsigma } \exp \left\{ -\frac{z^{2} }{2} \right\} dz\le
\\
\le \displaystyle\varepsilon L\frac{2}{\sqrt{2\pi } } \left[\int\limits_{0}^{1}z^{\varsigma } \exp \left\{ -\frac{z^{2} }{2} \right\} dz -\int\limits_{1}^{\infty }\exp \left\{ -\frac{z^{2} }{2} \right\}d\left( -\frac{z^{2} }{2}\right)\right]\le
\\
\le \displaystyle\varepsilon L\frac{2}{\sqrt{2\pi } } \Bigl[z^{\varsigma } |_{z=1} +1\Bigr]\le \varepsilon L\frac{4}{\sqrt{2\pi } }  .
\end{array}
\end{equation}
Так что при $\varepsilon \to 0$, $\varepsilon > 0$, для произвольных функций $f(t;x)$, которые для любых $t\in [0,T]$ ограничены с вероятностью 1 и   удовлетворяют условию Липшица по компоненте  $x$, справедливо равенство \eqref{GrindEQ__2_5_3_}.
\end{proof-p}

Дельта-функция и ее производные, примененные к ограниченной с вероятностью 1 функции, аналогичны свойствам, относящимся к детерминированным функциям. Поскольку применение формулы Ито к $\delta$-функции невозможно, в дальнейшем нам понадобится следующий результат.
\begin{proper}\label{lem-dif1}
Пусть функция $f(t;x)$ и ее производные первого и второго порядка по $x$  ограничены с вероятностью 1 для любых $t\in [0,T]$  и  удовлетворяют условию Липшица по компоненте  $x$. Тогда имеют место представления

\begin{equation}\label{a453a}
\begin{array}{c}
\displaystyle\frac{\partial }{\partial x} f(t;x)=-\displaystyle\lim\limits_{\varepsilon \downarrow 0} \int\limits_{-\infty }^{\infty } f(t;y)\frac{\partial }{\partial y}\delta_{\varepsilon}(y-x)  dy=\\
=\displaystyle\lim\limits_{\varepsilon \downarrow 0} \frac{1}{\varepsilon \sqrt{2\pi } } \int\limits_{-\infty }^{\infty }\exp \left\{ -\frac{(y-x)^{2} }{2\varepsilon ^{2} } \right\} \frac{\partial }{\partial y} f(t;y) dy,
\end{array}
\end{equation}
\begin{equation}\label{a453b}
\begin{array}{c} \displaystyle\frac{\partial^{2} }{\partial x^{2}} f(t;x)=\displaystyle\lim\limits_{\varepsilon \downarrow 0} \int\limits_{-\infty }^{\infty } f(t;y)\frac{\partial^{2} }{\partial y^{2}}\delta_{\varepsilon}(y-x)  dy=\\
=\displaystyle\lim\limits_{\varepsilon \downarrow 0} \frac{1}{\varepsilon \sqrt{2\pi } } \int\limits_{-\infty }^{\infty }\exp \left\{ -\frac{(y-x)^{2} }{2\varepsilon ^{2} } \right\}  \frac{\partial^{2} }{\partial y^{2}}f(t;y)  dy.
\end{array}
\end{equation}

\end{proper}

\begin{proof-p}
 Покажем, что и в этом случае справедливы оценки, подобные выполненным выше, на основе соответствующих предельных представлений. Продифференцируем \eqref{GrindEQ__2_5_3_}, применим интегрирование по частям:
$$
\begin{array}{c}
\displaystyle\frac{\partial }{\partial x} f(x)=\int\limits_{-\infty }^{\infty }f(y) \frac{\partial }{\partial x} \delta (y-x)dy=
 \\
  =\displaystyle\lim\limits_{\varepsilon \downarrow 0} \frac{1}{\varepsilon \sqrt{2\pi } } \int\limits_{-\infty }^{\infty }f(y) \frac{\partial }{\partial x} \exp \left\{ -\frac{(y-x)^{2} }{2\varepsilon ^{2} }\right\} dy=
  \\
  =-\displaystyle\lim\limits_{\varepsilon \downarrow 0} \frac{1}{\varepsilon \sqrt{2\pi } } \int\limits_{-\infty }^{\infty }f(y) \frac{\partial }{\partial y} \exp \left\{ -\frac{(y-x)^{2} }{2\varepsilon ^{2} }\right\} dy=
 \end{array}
$$
$$
\begin{array}{c}
 =-\displaystyle\lim\limits_{\varepsilon \downarrow 0} f(y)\frac{1}{\varepsilon \sqrt{2\pi } } \exp \left\{ -\frac{(y-x)^{2} }{2\varepsilon ^{2} } \right\} \Bigl. \Bigr|_{-\infty }^{+\infty } +\\
+\displaystyle\lim\limits_{\varepsilon \downarrow 0} \frac{1}{\varepsilon \sqrt{2\pi } } \int\limits_{-\infty }^{\infty } \exp \left\{ -\frac{(y-x)^{2} }{2\varepsilon ^{2} } \right\} \frac{\partial }{\partial y} f(y)dy,\ \ \ \varepsilon >0.
\end{array}
$$
При условии, что производная $f^{'}_{y} (y)$ удовлетворяет условию Липшица,  для асимптотического представления производной от $\delta (x)$, доказательство совпадает с предшествующим.
Аналогично устанавливается справедливость соотношений, где присутствуют вторые производные $\delta (x)$, и выполнено условие Липшица  для второй производной $f^{''}_{y} (y)$. Обобщая использовавшиеся ограничения видим, что в дальнейшем необходимо потребовать  существование  и непрерывность $f^{'}_{y} (y)$ и липшицевость для вторых производных.
\end{proof-p}

Перейдем к доказательству соотношения \eqref{GrindEQ__2_5_2_} теоремы, воспользовавшись установленным в утверждениях \ref{lem-eps} и \ref{lem-dif1} свойствами  $\delta (x)$ на классе всех непрерывных и ограниченных с вероятностью единица функций.

\begin{proofs}
Докажем, что при условии совпадения начальных значений
\begin{equation}\label{r01}
F(t;{\bf x}(t))= {\underset{\varepsilon \downarrow 0}{\mbox { \textrm{l.i.m.}}}} \,
F_{\varepsilon}(t;{\bf x}(t)),
\end{equation}
где
\begin{equation} \label{GrindEQ__2_5_4_}
F_{\varepsilon}(t;{\bf x}(t))=\displaystyle\int\limits_{\mathbb{R}^{n} }\prod_{i=1}^{n}\delta_{\varepsilon} (y_{i} -x_{i} (t))  F(t,{\bf y})d\Gamma ({\bf y}),
\end{equation}
и $d\Gamma({\bf y})=\displaystyle\prod\limits_{i=1}^{n}dy_{i}$  элемент фазового объема,  а $F(t;{\bf x}(t))$ -- решение \eqref{GrindEQ__2_5_2_}.

Утверждение \eqref{r01} равносильно  равенству
\begin{equation} \label{M-ots}
\lim\limits_{\varepsilon\downarrow 0}\mathbf{M}\Bigl[
\bigl|F_{\varepsilon}(t;{\bf x}(t))-F(t;{\bf x}(t))\bigr|^{2}
\Bigr]=0,
\end{equation}
которое будем доказывать.

Воспользовавшись обобщенной формулой Ито \eqref{Ayd2}, продифференцируем выражение $\displaystyle\prod\limits_{i=1}^{n}\delta_{\varepsilon}(y_{i}-x_{i}(t))F(t,{\bf y})$ в \eqref{GrindEQ__2_5_4_}. Положим: $h(t;\xi)=h(t;{\bf x};F)$ и $h(t;\xi(t))=h\bigl(t;{\bf x}(t);F(t;{\bf y})\bigr)=\displaystyle\prod\limits_{i=1}^{n}\delta_{\varepsilon}(y_{i}-x_{i}(t))F(t,{\bf y})$, где ${\bf x}(t)$ -- подчинено уравнению \eqref{GrindEQ__2_1_1_}, а
$F(t;{\bf x})$ -- уравнению \eqref{GrindEQ__2_5_1_}.
Вычислим последовательно все слагаемые для $d_{t}h(t;\xi(t))$. Вернемся для удобства к использованию знака суммы.

{\bf 1.} Вычислим слагаемое $h_{t}^{'}(t;\xi(t))$. В силу того, что в этой функции нет явной зависимости от $t$:
\begin{equation}\label{sl-1}
h_{t}^{'}(t;\xi(t))= \dfrac{\partial}{\partial t}\bigl[\prod\limits_{i=1}^{n}\delta_{\varepsilon}(y_{i}-x_{i})\cdot F\bigr]\bigl|_{\xi=\xi(t)}=0.
\end{equation}
 Здесь и далее в этой части условие $\xi=\xi(t)$ означает, подстановку $x_{i}=x_{i}(t)$, $F=F(t;{\bf y})$.

{\bf 2.} Вычислим сумму $\displaystyle\sum\limits_{i}h_{\xi_{i}}^{'}(t;\xi(t))\cdot \widetilde{a}_{i}(t)$:
\begin{equation}\label{sl-2}
\begin{array}{c}
\displaystyle\sum\limits_{i}h_{\xi_{i}}^{'}(t;\xi(t))\cdot \widetilde{a}_{i}(t)=
\sum\limits_{i}^{n}{a}_{i}(t)\cdot F\cdot \dfrac{\partial}{\partial x_{i}}\prod\limits_{i=1}^{n}\delta_{\varepsilon}(y_{i}-x_{i})\bigl|_{\xi=\xi(t)} +\\
+Q(t,{\bf y})\displaystyle\prod\limits_{i=1}^{n}\delta_{\varepsilon}(y_{i}-x_{i})\dfrac{\partial F}{\partial F}\bigl|_{\xi=\xi(t)}=\\
=\displaystyle\sum\limits_{i}^{n}{a}_{i}(t)F(t,{\bf y})\dfrac{\partial}{\partial x_{i}}\prod\limits_{i=1}^{n}\delta_{\varepsilon}(y_{i}-x_{i})\bigl|_{\xi=\xi(t)} +Q(t,{\bf y})\displaystyle\prod\limits_{i=1}^{n}\delta_{\varepsilon}(y_{i}-x_{i}(t)).
\end{array}
\end{equation}

{\bf 3.} Вычислим сумму $\displaystyle \sum\limits_{i=1}
 \sum\limits_{j=1}
 \sum\limits_{k=1}^{m}\widetilde{b}_{i\,k}(t;\xi(t))\widetilde{b}_{j\,k}(t;\xi(t))
 \frac{\partial^{\,2} h(t;\xi) }{\partial \xi_{i}\partial\xi_{j}}\bigl|_{\xi =\xi(t)}$:
\begin{equation}\label{sl-3}
\begin{array}{c}
\displaystyle \sum\limits_{i=1}
 \sum\limits_{j=1}
 \sum\limits_{k=1}^{m}\widetilde{b}_{i\,k}(t;\xi(t))\widetilde{b}_{j\,k}(t;\xi(t))
\frac{\partial^{\,2} h(t;\xi) }{\partial \xi_{i}\partial\xi_{j}}\bigl|_{\xi =\xi(t)}=\\
  =\displaystyle \sum\limits_{i=1}^{n}
 \sum\limits_{j=1}^{n}
 \sum\limits_{k=1}^{m}{b}_{i\,k}(t){b}_{j\,k}(t) F(t,{\bf y})
 \frac{\partial^{\,2}}{\partial x_{i}\partial x_{j}}\prod\limits_{i=1}^{n}\delta_{\varepsilon}(y_{i}-x_{i})\bigl|_{\xi =\xi(t)} +\\
 + 2\cdot\displaystyle\sum\limits_{i=1}
 \sum\limits_{k=1}^{m}{b}_{i\,k}(t)D_{k}(t,{\bf y})
\frac{\partial}{\partial x_{i}}\prod\limits_{i=1}^{n}\delta_{\varepsilon}(y_{i}-x_{i}) \bigl|_{\xi =\xi(t)}.
\end{array}
\end{equation}

{\bf 4.} Вычислим винеровскую составляющую.
\begin{equation}\label{sl-4}
\begin{array}{c}
\displaystyle \sum\limits_{i=1}
  \sum\limits_{k=1}^{m}h^{'}_{\xi_{i}}(t;\xi)\widetilde{b}_{i\,k}(t;\xi(t))=\\
  =F(t;{\bf y})\displaystyle \sum\limits_{i=1}
  \sum\limits_{k=1}^{m}{b}_{i\,k}(t)\frac{\partial}{\partial x_{i}}\prod\limits_{i=1}^{n}\delta_{\varepsilon}(y_{i}-x_{i}) \bigl|_{\xi =\xi(t)}+\\
  +\displaystyle\sum\limits_{k=1}^{m}D_{k}(t;{\bf y})\prod\limits_{i=1}^{n}\delta_{\varepsilon}(y_{i}-x_{i}(t)).
\end{array}
\end{equation}

{\bf 5.} Вычислим пуассоновскую составляющую. С учетом, что
$$
h(t;\xi(t)+\pi(t))=h\bigl(t;{\bf x}(t)+\pi_{1}(t);F(t;{\bf y})+\pi_{2}(t)\bigr),
$$
получаем:
\begin{equation}\label{sl-5}
\begin{array}{c}
\displaystyle\int\Bigl[h(t;\xi(t)+\pi(t))-h(t;\xi(t))\Bigr]\nu(dt;d\gamma)=\\
=\displaystyle\int\Bigl[
\prod\limits_{i=1}^{n}\delta_{\varepsilon}\bigl(y_{i}-x_{i}(t)-g_{i}(t;\gamma)\bigr)
\bigl(F(t,{\bf y})+G(t,{\bf y};\gamma)\bigr)-\\
-
\displaystyle\prod\limits_{i=1}^{n}\delta_{\varepsilon}(y_{i}-x_{i}(t))F(t,{\bf y})\Bigr]\nu(dt;d\gamma).
\end{array}
\end{equation}

В итоге, с учетом полученных выражений \eqref{sl-1}-\eqref{sl-5}, приходим к следующему результату (знаки суммирования  опустим, и в дальнейшем будем иметь в
виду, что суммирование проводится по индексам, встречающимся дважды):
\begin{equation}\label{itog-eps}
\begin{array}{c}
d_{t}\Bigl[\displaystyle\prod\limits_{i=1}^{n}\delta_{\varepsilon}(y_{i}-x_{i}(t))F(t,{\bf y})\Bigr]=\\
=\Bigl[{a}_{i}(t)F(t,{\bf y})\displaystyle\dfrac{\partial}{\partial x_{i}}\prod\limits_{i=1}^{n}\delta_{\varepsilon}(y_{i}-x_{i})\bigl|_{\xi=\xi(t)} +Q(t,{\bf y})\displaystyle\prod\limits_{i=1}^{n}\delta_{\varepsilon}(y_{i}-x_{i}(t))+\\
+ \displaystyle\frac{1}{2}
{b}_{i\,k}(t){b}_{j\,k}(t) F(t,{\bf y})
 \frac{\partial^{\,2}}{\partial x_{i}\partial x_{j}}\prod\limits_{i=1}^{n}\delta_{\varepsilon}(y_{i}-x_{i})\bigl|_{\xi =\xi(t)} +\\
 +\displaystyle{b}_{i\,k}(t)D_{k}(t,{\bf y})
\frac{\partial}{\partial x_{i}}\prod\limits_{i=1}^{n}\delta_{\varepsilon}(y_{i}-x_{i}) \bigl|_{\xi =\xi(t)}\Bigr]dt+\\
+\Bigl[F(t;{\bf y})\displaystyle {b}_{i\,k}(t)\frac{\partial}{\partial x_{i}}\prod\limits_{i=1}^{n}\delta_{\varepsilon}(y_{i}-x_{i}) \bigl|_{\xi =\xi(t)}+\\
  +\displaystyle D_{k}(t;{\bf y})\prod\limits_{i=1}^{n}\delta_{\varepsilon}(y_{i}-x_{i}(t))\Bigr]dw_{k}(t)+\\
+\displaystyle\int\Bigl[
\prod\limits_{i=1}^{n}\delta_{\varepsilon}\bigl(y_{i}-x_{i}(t)-g_{i}(t;\gamma)\bigr)
\bigl(F(t,{\bf y})+G(t,{\bf y};\gamma)\bigr)-\\
-
\displaystyle\prod\limits_{i=1}^{n}\delta_{\varepsilon}(y_{i}-x_{i}(t))F(t,{\bf y})\Bigr]\nu(dt;d\gamma).
\end{array}
\end{equation}

Перейдем в правой части равенства от частных производных по компонентам вектора ${\bf x}$ к частным производных по компонентам  вектора ${\bf y}$.
\begin{equation}\label{y-eps}
\begin{array}{c}
d_{t}\Bigl[\displaystyle\prod\limits_{i=1}^{n}\delta_{\varepsilon}(y_{i}-x_{i}(t))F(t,{\bf y})\Bigr]=\\
=\Bigl[-{a}_{i}(t)F(t,{\bf y})\displaystyle\dfrac{\partial}{\partial y_{i}}\prod\limits_{i=1}^{n}\delta_{\varepsilon}(y_{i}-x_{i}(t) +Q(t,{\bf y})\displaystyle\prod\limits_{i=1}^{n}\delta_{\varepsilon}(y_{i}-x_{i}(t))+\\
+ \displaystyle\frac{1}{2}
{b}_{i\,k}(t){b}_{j\,k}(t) F(t,{\bf y})
 \frac{\partial^{\,2}}{\partial y_{i}\partial y_{j}}\prod\limits_{i=1}^{n}\delta_{\varepsilon}(y_{i}-x_{i}(t)) -\\
 -\displaystyle{b}_{i\,k}(t)D_{k}(t,{\bf y})
\frac{\partial}{\partial y_{i}}\prod\limits_{i=1}^{n}\delta_{\varepsilon}(y_{i}-x_{i}(t)) \Bigr]dt+\\
+\Bigl[-F(t;{\bf y})\displaystyle {b}_{i\,k}(t)\frac{\partial}{\partial y_{i}}\prod\limits_{i=1}^{n}\delta_{\varepsilon}(y_{i}-x_{i}(t)) +\\
  +\displaystyle D_{k}(t;{\bf y})\prod\limits_{i=1}^{n}\delta_{\varepsilon}(y_{i}-x_{i}(t))\Bigr]dw_{k}(t)+\\
+\displaystyle\int\Bigl[
\prod\limits_{i=1}^{n}\delta_{\varepsilon}\bigl(y_{i}-x_{i}(t)-g_{i}(t;\gamma)\bigr)
\bigl(F(t,{\bf y})+G(t,{\bf y};\gamma)\bigr)-\\
-
\displaystyle\prod\limits_{i=1}^{n}\delta_{\varepsilon}(y_{i}-x_{i}(t))F(t,{\bf y})\Bigr]\nu(dt;d\gamma).
\end{array}
\end{equation}

Полученное ОСДУ \eqref{y-eps} запишем в интегральном виде.
и возвратимся к интегрированию по пространству $\mathbb{R}^{n} $, с учетом  равенства~\eqref{GrindEQ__2_5_4_}.
\begin{equation}\label{integr-prostr}
\begin{array}{c}
F_{\varepsilon}(t;{\bf x}(t))=\displaystyle\int\limits_{\mathbb{R}^{n} }\prod\limits_{i=1}^{n}\delta_{\varepsilon}(y_{i}-x_{i}(t))F(t,{\bf y})d\Gamma({\bf y})=\\
=\displaystyle\int\limits_{\mathbb{R}^{n} }\prod\limits_{i=1}^{n}\delta_{\varepsilon}(y_{i}^{0}-x_{i}(0))F(0,{\bf y_{0}})d\Gamma({\bf y})-\\
-\displaystyle\int\limits_{\mathbb{R}^{n} }\displaystyle\int\limits_{0}^{t}{a}_{i}(\tau)F(\tau,{\bf y})\dfrac{\partial}{\partial y_{i}}\prod\limits_{i=1}^{n}\delta_{\varepsilon}(y_{i}-x_{i}(\tau))d\tau d\Gamma({\bf y}) +\\
+\displaystyle\int\limits_{\mathbb{R}^{n} }\displaystyle\int\limits_{0}^{t}Q(\tau,{\bf y})\displaystyle\prod\limits_{i=1}^{n}\delta_{\varepsilon}(y_{i}-x_{i}(\tau))d\tau d\Gamma({\bf y})+\\
+ \displaystyle\frac{1}{2}\displaystyle\int\limits_{\mathbb{R}^{n}}
\int\limits_{0}^{t}{b}_{i\,k}(\tau){b}_{j\,k}(\tau) F(\tau,{\bf y})
 \frac{\partial^{\,2}}{\partial y_{i}\partial y_{j}}\prod\limits_{i=1}^{n}\delta_{\varepsilon}(y_{i}-x_{i}(\tau))d\tau d\Gamma({\bf y}) -\\
 -\displaystyle\int\limits_{\mathbb{R}^{n}}\displaystyle\int\limits_{0}^{t}{b}_{i\,k}(\tau)D_{k}(\tau,{\bf y})
\frac{\partial}{\partial y_{i}}\prod\limits_{i=1}^{n}\delta_{\varepsilon}(y_{i}-x_{i}(\tau)) d\tau d\Gamma({\bf y})-
\end{array}
\end{equation}
$$
\begin{array}{c}
-\displaystyle\int\limits_{\mathbb{R}^{n}}\int\limits_{0}^{t}F(\tau;{\bf y})\displaystyle {b}_{i\,k}(\tau)\frac{\partial}{\partial y_{i}}\prod\limits_{i=1}^{n}\delta_{\varepsilon}(y_{i}-x_{i}(\tau)) dw_{k}(\tau)d\Gamma({\bf y})+\\
  +\displaystyle\int\limits_{\mathbb{R}^{n}}\displaystyle\int\limits_{0}^{t}D_{k}(\tau;{\bf y})\prod\limits_{i=1}^{n}\delta_{\varepsilon}(y_{i}-x_{i}(\tau))dw_{k}(\tau)d\Gamma({\bf y})+
\\
+\displaystyle\int\limits_{\mathbb{R}^{n}}\int\limits_{0}^{t}\int\Bigl[
\prod\limits_{i=1}^{n}\delta_{\varepsilon}\bigl(y_{i}-x_{i}(\tau)-g_{i}(\tau;\gamma)\bigr)
\bigl(F(\tau,{\bf y})+G(\tau,{\bf y};\gamma)\bigr)-\\
-
\displaystyle\prod\limits_{i=1}^{n}\delta_{\varepsilon}(y_{i}-x_{i}(\tau))F(\tau,{\bf y})\Bigr]\nu(d\tau;d\gamma)d\Gamma({\bf y}).
\end{array}
$$
Запишем ОСДУ \eqref{GrindEQ__2_5_2_} в интегральном виде.
\begin{equation}\label{d2}
\begin{array}{c}
F(t;{\bf x}(t))=F(0;{\bf x}(0))+\displaystyle\int\limits_{0}^{t}Q(\tau;{\bf x}(\tau))d\tau
+\int\limits_{0}^{t}D_{k} (\tau;{\bf x}(\tau))dw_{k}(\tau) +\\
+\displaystyle\int\limits_{0}^{t}b_{i,k} (\tau)\displaystyle\frac{\partial F(\tau;{\bf x})}{\partial x_{i} }\bigl|_{{\bf x}={\bf x}(\tau) }dw_{k}(\tau) +
\\
+\displaystyle \int\limits_{0}^{t}\Bigl[a_{i} (\tau)\frac{\partial F(\tau;{\bf x}) }{\partial x_{i} }\bigl|_{{\bf x}={\bf x}(\tau) } +\frac{1}{2} b_{i,k} (\tau)b_{j,k} (\tau)\frac{\partial ^{2} F(\tau;{\bf x})}{\partial x_{i} \partial x_{j} }\bigl|_{{\bf x}={\bf x}(\tau) }\Bigr. + \\
+\Bigl.\displaystyle
b_{i,k} (\tau)\frac{\partial D_{k} (\tau;{\bf x})}{\partial x_{i} }\bigl|_{{\bf x}={\bf x}(\tau) } \Bigr]d\tau+ \\
+\displaystyle\int\limits_{0}^{t}\int \Bigl[F\bigl(\tau;{\bf x}(\tau)+g(\tau;\gamma )\bigr)-F(\tau;{\bf x}(\tau))\Bigr]\nu (d\tau;d\gamma ) +\\
+ \displaystyle\int\limits_{0}^{t}\int G\bigl(\tau;{\bf x}(\tau)+g(\tau;\gamma );\gamma\bigr)\nu (d\tau;d\gamma ).
\end{array}
\end{equation}

Рассмотрим разность выражений \eqref{integr-prostr} и \eqref{d2} с учетом допредельных свойств $\delta $-функции и возможности изменения порядка интегрирования:
\begin{equation}\label{integr-prostr00}
\begin{array}{c}
F_{\varepsilon}(t;{\bf x}(t))-F(t;{\bf x}(t))=
F_{\varepsilon}(0,{\bf x}(0))-F(0;{\bf x}(0))+\\
+\displaystyle\int\limits_{0}^{t}a_{i} (\tau)\Bigl[-\displaystyle\int\limits_{\mathbb{R}^{n}}F(\tau,{\bf y})\dfrac{\partial}{\partial y_{i}}\prod\limits_{i=1}^{n}\delta_{\varepsilon}(y_{i}-x_{i}(\tau)) d\Gamma({\bf y})-\displaystyle \frac{\partial F(\tau;{\bf x}) }{\partial x_{i} }\bigl|_{{\bf x}={\bf x}(\tau) }\Bigr]d\tau+\\
+\displaystyle\int\limits_{0}^{t}\Bigl[\displaystyle\int\limits_{\mathbb{R}^{n}}Q(\tau,{\bf y})\displaystyle\prod\limits_{i=1}^{n}\delta_{\varepsilon}(y_{i}-x_{i}(\tau)) d\Gamma({\bf y})-Q(\tau;{\bf x}(\tau))\Bigr]d\tau+
\end{array}
\end{equation}
\begin{equation*}
\begin{array}{c}
+ \displaystyle\frac{1}{2}\displaystyle
\int\limits_{0}^{t}{b}_{i\,k}(\tau){b}_{j\,k}(\tau)\Bigl[\int\limits_{\mathbb{R}^{n}} F(\tau,{\bf y})
 \frac{\partial^{\,2}}{\partial y_{i}\partial y_{j}}\prod\limits_{i=1}^{n}\delta_{\varepsilon}(y_{i}-x_{i}(\tau)) d\Gamma({\bf y}) -\\
-\displaystyle\frac{\partial ^{2} F(\tau;{\bf x})}{\partial x_{i} \partial x_{j} }\bigl|_{{\bf x}={\bf x}(\tau) }\Bigr]d\tau +\\
  +\displaystyle\int\limits_{0}^{t}{b}_{i\,k}(\tau)\Bigl[-\int\limits_{\mathbb{R}^{n}}D_{k}(\tau,{\bf y})
\frac{\partial}{\partial y_{i}}\prod\limits_{i=1}^{n}\delta_{\varepsilon}(y_{i}-x_{i}(\tau)) d\Gamma({\bf y})-\\
-\displaystyle\frac{\partial D_{k} (\tau;{\bf x})}{\partial x_{i} }\bigl|_{{\bf x}={\bf x}(\tau) }\Bigr]d\tau+\\
-\displaystyle\int\limits_{0}^{t}{b}_{i\,k}(\tau)\Bigl[-\int\limits_{\mathbb{R}^{n}}F(\tau;{\bf y})\displaystyle \frac{\partial}{\partial y_{i}}\prod\limits_{i=1}^{n}\delta_{\varepsilon}(y_{i}-x_{i}(\tau)) d\Gamma({\bf y})-\\
-\displaystyle\frac{\partial F(\tau;{\bf x})}{\partial x_{i} }\bigl|_{{\bf x}={\bf x}(\tau) }\Bigr]dw_{k}(\tau) +
\end{array}
\end{equation*}
\begin{equation*}
\begin{array}{c}
  +\displaystyle\int\limits_{0}^{t}\Bigl[\int\limits_{\mathbb{R}^{n}}D_{k}(\tau;{\bf y})\prod\limits_{i=1}^{n}\delta_{\varepsilon}(y_{i}-x_{i}(\tau))d\Gamma({\bf y})
-\displaystyle D_{k} (\tau;{\bf x}(\tau))\Bigr]dw_{k}(\tau) +\\
+\displaystyle\int\limits_{0}^{t}\int\Bigl(\int\limits_{\mathbb{R}^{n}}\Bigl[
\prod\limits_{i=1}^{n}\delta_{\varepsilon}\bigl(y_{i}-x_{i}(\tau)-g_{i}(\tau;\gamma)\bigr)
\bigl(F(\tau,{\bf y})+G(\tau,{\bf y};\gamma)\bigr)-\\
-
\displaystyle\prod\limits_{i=1}^{n}\delta_{\varepsilon}(y_{i}-x_{i}(\tau))F(\tau,{\bf y})\Bigr]d\Gamma({\bf y})-\\
- F\bigl(\tau;{\bf x}(\tau)+g(\tau;\gamma )\bigr)+F(\tau;{\bf x}(\tau))
- G\bigl(\tau;{\bf x}(\tau)+g(\tau;\gamma );\gamma\bigr)\Bigr)\nu (d\tau;d\gamma ).
\end{array}
\end{equation*}

Поскольку дифференцирование по Ито дельта-функции невозможно, воспользуемся результатами микролемм -- равенствами \eqref{GrindEQ__2_5_3_}, \eqref{a453a} и \eqref{a453b}:
\begin{equation}\label{integr-prostr11}
\begin{array}{c}
\Bigl|F_{\varepsilon}(t;{\bf x}(t))-F(t;{\bf x}(t))\Bigr|\leq
\Bigl|F_{\varepsilon}(0,{\bf x}(0))-F(0;{\bf x}(0))\Bigr|+\\
+\Biggl|\displaystyle\int\limits_{0}^{t}a_{i} (\tau)\Bigl[\dfrac{1}{(\varepsilon\sqrt{2\pi})^{n}}\displaystyle\int\limits_{\mathbb{R}^{n} }\prod\limits_{i=1}^{n}\exp\left\{-\frac{(y_{i}-x_{i}(\tau))^{2}}{2\varepsilon^{2}}   \right\} \dfrac{\partial}{\partial y_{i}}F(\tau,{\bf y}) d\Gamma({\bf y})-\\
-\displaystyle \frac{\partial F(\tau;{\bf x}) }{\partial x_{i} }\bigl|_{{\bf x}={\bf x}(\tau) }\Bigr]d\tau\Biggr|+\\
+\Biggl|\displaystyle\int\limits_{0}^{t}\Bigl[\dfrac{1}{(\varepsilon\sqrt{2\pi})^{n}}
\displaystyle\int\limits_{\mathbb{R}^{n}}\prod\limits_{i=1}^{n}\exp\left\{-\frac{(y_{i}-x_{i}(\tau))^{2}}{2\varepsilon^{2}}   \right\} Q(\tau,{\bf y}) d\Gamma({\bf y})-\\
-Q(\tau;{\bf x}(\tau))\Bigr]d\tau\Biggr|+
\end{array}
\end{equation}
\begin{equation*}
\begin{array}{c}
+\Biggl| \displaystyle\frac{1}{2}\displaystyle
\int\limits_{0}^{t}\Bigl[\dfrac{1}{(\varepsilon\sqrt{2\pi})^{n}}\displaystyle\int\limits_{\mathbb{R}^{n} }\prod\limits_{i=1}^{n}\exp\left\{-\frac{(y_{i}-x_{i}(\tau))^{2}}{2\varepsilon^{2}}   \right\}
 \frac{\partial^{\,2}}{\partial y_{i}\partial y_{j}}F(\tau,{\bf y}) d\Gamma({\bf y}) -\\
-\displaystyle\frac{\partial ^{2} F(\tau;{\bf x})}{\partial x_{i} \partial x_{j} }\bigl|_{{\bf x}={\bf x}(\tau) }\Bigr]{b}_{i\,k}(\tau){b}_{j\,k}(\tau)d\tau \Biggr|+
\end{array}
\end{equation*}
\begin{equation*}
\begin{array}{c}
  +\Biggl|\displaystyle\int\limits_{0}^{t}{b}_{i\,k}(\tau)\Bigl[
  \dfrac{1}{(\varepsilon\sqrt{2\pi})^{n}}\displaystyle\int\limits_{\mathbb{R}^{n} }\prod\limits_{i=1}^{n}\exp\left\{-\frac{(y_{i}-x_{i}(\tau))^{2}}{2\varepsilon^{2}}   \right\} \frac{\partial}{\partial y_{i}}D_{k}(\tau,{\bf y}) d\Gamma({\bf y})-\\
-\displaystyle\frac{\partial D_{k} (\tau;{\bf x})}{\partial x_{i} }\bigl|_{{\bf x}={\bf x}(\tau) }\Bigr]d\tau\Biggr|+
\end{array}
\end{equation*}
\begin{equation*}
\begin{array}{c}
+\Biggl|\displaystyle\int\limits_{0}^{t}{b}_{i\,k}(\tau)\Bigl[
\dfrac{1}{(\varepsilon\sqrt{2\pi})^{n}}\displaystyle\int\limits_{\mathbb{R}^{n} }\prod\limits_{i=1}^{n}\exp\left\{-\frac{(y_{i}-x_{i}(\tau))^{2}}{2\varepsilon^{2}}   \right\}\frac{\partial}{\partial y_{i}}F(\tau;{\bf y}) d\Gamma({\bf y})-\\
-\displaystyle\frac{\partial F(\tau;{\bf x})}{\partial x_{i} }\bigl|_{{\bf x}={\bf x}(\tau) }\Bigr]dw_{k}(\tau)\Biggr| +\\
  +\Biggl|\displaystyle\int\limits_{0}^{t}\Bigl[
  \dfrac{1}{(\varepsilon\sqrt{2\pi})^{n}}\displaystyle\int\limits_{\mathbb{R}^{n} }\prod\limits_{i=1}^{n}\exp\left\{-\frac{(y_{i}-x_{i}(\tau))^{2}}{2\varepsilon^{2}}   \right\}D_{k}(\tau;{\bf y})d\Gamma({\bf y})-\\
-\displaystyle D_{k} (\tau;{\bf x}(\tau))\Bigr]dw_{k}(\tau)\Biggr| +
\end{array}
\end{equation*}
\begin{equation*}
\begin{array}{c}
+\Biggl|\displaystyle\int\limits_{0}^{t}\int\Bigl[
\dfrac{1}{(\varepsilon\sqrt{2\pi})^{n}}\displaystyle\int\limits_{\mathbb{R}^{n} }\prod\limits_{i=1}^{n}\exp\left\{-\frac{\bigl(y_{i}-x_{i}(\tau)-g_{i}(\tau;\gamma)\bigr)^{2}}{2\varepsilon^{2}}   \right\}\cdot\\
\cdot F(\tau,{\bf y})d\Gamma({\bf y})  -F\bigl(\tau;{\bf x}(\tau)+g(\tau;\gamma )\bigr)\Bigr]\nu (d\tau;d\gamma )\Biggr|+\\
+
\Biggl|\displaystyle\int\limits_{0}^{t}\int\Bigl[
\dfrac{1}{(\varepsilon\sqrt{2\pi})^{n}}\displaystyle\int\limits_{\mathbb{R}^{n} }\prod\limits_{i=1}^{n}\exp\left\{-\frac{(y_{i}-x_{i}(\tau))^{2}}{2\varepsilon^{2}}   \right\} F(\tau,{\bf y})d\Gamma({\bf y})-\\
-F(\tau;{\bf x}(\tau))\Bigr]\nu (d\tau;d\gamma )\Biggr|+ \\
+
\Biggl|\displaystyle\int\limits_{0}^{t}\int\Bigl[
\dfrac{1}{(\varepsilon\sqrt{2\pi})^{n}}\displaystyle\int\limits_{\mathbb{R}^{n} }\prod\limits_{i=1}^{n}\exp\left\{-\frac{\bigl(y_{i}-x_{i}(\tau)-g_{i}(\tau;\gamma)\bigr)^{2}}{2\varepsilon^{2}}   \right\}\cdot\\
\cdot
G(\tau,{\bf y};\gamma)\bigr)d\Gamma({\bf y})
-G\bigl(\tau;{\bf x}(\tau)+g(\tau;\gamma );\gamma\bigr)\Bigr]\nu (d\tau;d\gamma )\Biggr|.
\end{array}
\end{equation*}

Интегралы, содержащие экспоненту, проинтегрируем по частям, проведем действия, аналогичные \eqref{dob1}, и произведем оценки сверху, аналогичные \eqref{ocenka1}.

\begin{equation}\label{ots-Q}
\begin{array}{c}
\Biggl|\displaystyle\int\limits_{0}^{t}\Bigl[\dfrac{1}{(\varepsilon\sqrt{2\pi})^{n}}\displaystyle
\int\limits_{\mathbb{R}^{n}}\prod\limits_{i=1}^{n}\exp\left\{-\frac{(y_{i}-x_{i}(\tau))^{2}}{2\varepsilon^{2}}   \right\} Q(\tau,{\bf y}) d\Gamma({\bf y})-\\
-Q(\tau;{\bf x}(\tau))\Bigr]d\tau\Biggr|\leq\\
\leq
\displaystyle\int\limits_{0}^{t}\Biggl|\dfrac{1}{(\sqrt{2\pi})^{n}}\displaystyle\int\limits_{\mathbb{R}^{n} }\prod\limits_{i=1}^{n}\exp\left\{-\frac{z_{i}^{2}(\tau)}{2}   \right\}\Bigl[ Q(\tau,\varepsilon  {\bf z}(\tau)+ {\bf x}(\tau))-\\
-Q(\tau;{\bf x}(\tau))\Bigr]  d\Gamma({\bf z})\Biggr|d\tau\leq\\
\leq
\displaystyle\int\limits_{0}^{t}\dfrac{\varepsilon^{n}2^{n+1}}{(\sqrt{2\pi})^{n}}L^{(1)}(\tau)d\tau= \dfrac{\varepsilon^{n}2^{n+1}}{(\sqrt{2\pi})^{n}}\displaystyle\int\limits_{0}^{t}L^{(1)}(\tau)d\tau.
\end{array}
\end{equation}
Аналогично:
\begin{equation}\label{ots-D}
\begin{array}{c}
\Biggl|\displaystyle\int\limits_{0}^{t}\Bigl[\dfrac{1}{(\varepsilon\sqrt{2\pi})^{n}}\displaystyle
\int\limits_{\mathbb{R}^{n} }\prod\limits_{i=1}^{n}\exp\left\{-\frac{(y_{i}-x_{i}(\tau))^{2}}{2\varepsilon^{2}}   \right\} D_{k}(\tau,{\bf y}) d\Gamma({\bf y})-\\
-D_{k}(\tau;{\bf x}(\tau))\Bigr]dw_{k}(\tau)\Biggr|\leq\\
\leq
\displaystyle\int\limits_{0}^{t}\dfrac{\varepsilon^{n}2^{n+1}}{(\sqrt{2\pi})^{n}}L^{(2)}(\tau)dw_{k}(\tau)= \dfrac{\varepsilon^{n}2^{n+1}}{(\sqrt{2\pi})^{n}}\displaystyle\int\limits_{0}^{t}L^{(2)}(\tau)dw_{k}(\tau).
\end{array}
\end{equation}
Далее:
\begin{equation}\label{ots-Fg}
\begin{array}{c}
\Biggl|\displaystyle\int\limits_{0}^{t}\int\Bigl[\dfrac{1}{(\varepsilon\sqrt{2\pi})^{n}}
\displaystyle\int\limits_{\mathbb{R}^{n}}\prod\limits_{i=1}^{n}\exp\left\{-\frac{\bigl(y_{i}-x_{i}(\tau)-g_{i}(\tau;\gamma)\bigr)^{2}}{2\varepsilon^{2}}   \right\}\cdot\\
\cdot F(\tau,{\bf y})d\Gamma({\bf y})  -F\bigl(\tau;{\bf x}(\tau)+g(\tau;\gamma )\bigr)\Bigr]\nu (d\tau;d\gamma )\Biggr|\leq \\
\leq
\dfrac{\varepsilon^{n}2^{n+1}}{(\sqrt{2\pi})^{n}}\displaystyle\int\limits_{0}^{t}\int
L^{(3)}(\tau;\gamma)\nu (d\tau;d\gamma ),
\end{array}
\end{equation}
\begin{equation}\label{ots-G}
\begin{array}{c}
\Biggl|\displaystyle\int\limits_{0}^{t}\int\Bigl[\dfrac{1}{(\varepsilon\sqrt{2\pi})^{n}}
\displaystyle\int\limits_{\mathbb{R}^{n} }\prod\limits_{i=1}^{n}\exp\left\{-\frac{\bigl(y_{i}-x_{i}(\tau)-g_{i}(\tau;\gamma)\bigr)^{2}}{2\varepsilon^{2}}   \right\}\cdot\\
\cdot
G(\tau,{\bf y};\gamma)\bigr)d\Gamma({\bf y})
-G\bigl(\tau;{\bf x}(\tau)+g(\tau;\gamma );\gamma\bigr)\Bigr]\nu (d\tau;d\gamma )\Biggr|\leq\\
\leq
\dfrac{\varepsilon^{n}2^{n+1}}{(\sqrt{2\pi})^{n}}\displaystyle\int\limits_{0}^{t}\int
L^{(4)}(\tau;\gamma)\nu (d\tau;d\gamma ).
\end{array}
\end{equation}
Учитывая замену $z=\dfrac{y-x}{\varepsilon}$, получаем, что
$$
\frac{\partial f(\varepsilon z + x)}{\partial z}=\varepsilon \, \frac{\partial f(\varepsilon z + x)}{\partial x},\ \ \ \frac{\partial^{2} f(\varepsilon z + x)}{\partial z^{2}}=\varepsilon^{2} \, \frac{\partial^{2} f(\varepsilon z + x)}{\partial x^{2}}.
$$
Тогда для выражений с первыми производными проведем действия, аналогичные \eqref{dob1}, и  оценки, аналогичные \eqref{ocenka1}:
\begin{equation}\label{ots-der-a}
\begin{array}{c}
\Biggl|\displaystyle\int\limits_{0}^{t}a_{i} (\tau)\Bigl[\dfrac{1}{(\varepsilon\sqrt{2\pi})^{n}}\displaystyle\int\limits_{\mathbb{R}^{n} }\prod\limits_{i=1}^{n}\exp\left\{-\frac{(y_{i}-x_{i}(\tau))^{2}}{2\varepsilon^{2}}   \right\} \dfrac{\partial}{\partial y_{i}}F(\tau,{\bf y}) d\Gamma({\bf y})-\\
-\displaystyle \frac{\partial F(\tau;{\bf x}) }{\partial x_{i} }\bigl|_{{\bf x}={\bf x}(\tau) }\Bigr]d\tau\Biggr|
\leq
\dfrac{\varepsilon^{n}2^{n+1}}{(\sqrt{2\pi})^{n}}
\displaystyle\int\limits_{0}^{t}\Bigl|{a}_{i}(\tau)\Bigr|
L^{(5)}(\tau)d\tau,
\end{array}
\end{equation}
\begin{equation}\label{ots-der-D}
\begin{array}{c}
\Biggl|\displaystyle\int\limits_{0}^{t}{b}_{i\,k}(\tau)\Bigl[\dfrac{1}{(\varepsilon\sqrt{2\pi})^{n}}\displaystyle\int\limits_{\mathbb{R}^{n} }\prod\limits_{i=1}^{n}\exp\left\{-\frac{(y_{i}-x_{i}(\tau))^{2}}{2\varepsilon^{2}}   \right\} \frac{\partial}{\partial y_{i}}D_{k}(\tau,{\bf y}) d\Gamma({\bf y})-\\
-\displaystyle\frac{\partial D_{k} (\tau;{\bf x})}{\partial x_{i} }\bigl|_{{\bf x}={\bf x}(\tau) }\Bigr]d\tau\Biggr|
\leq
\dfrac{\varepsilon^{n}2^{n+1}}{(\sqrt{2\pi})^{n}}
\displaystyle\int\limits_{0}^{t}\Bigl|{b}_{i\,k}(\tau)\Bigr|
L^{(6)}(\tau)d\tau,
\end{array}
\end{equation}
\begin{equation}\label{ots-der-F}
\begin{array}{c}
\Biggl|\displaystyle\int\limits_{0}^{t}{b}_{i\,k}(\tau)\Bigl[\dfrac{1}{(\varepsilon\sqrt{2\pi})^{n}}\displaystyle\int\limits_{\mathbb{R}^{n} }\prod\limits_{i=1}^{n}\exp\left\{-\frac{(y_{i}-x_{i}(\tau))^{2}}{2\varepsilon^{2}}   \right\}\frac{\partial}{\partial y_{i}}F(\tau;{\bf y}) d\Gamma({\bf y})-\\
-\displaystyle\frac{\partial F(\tau;{\bf x})}{\partial x_{i} }\bigl|_{{\bf x}={\bf x}(\tau) }\Bigr]dw_{k}(\tau)\Biggr|
\leq
\dfrac{\varepsilon^{n}2^{n+1}}{(\sqrt{2\pi})^{n}}
\displaystyle\int\limits_{0}^{t}\Bigl|{b}_{i\,k}(\tau)\Bigr|
L^{(7)}(\tau)dw_{k}(\tau).
\end{array}
\end{equation}
Для выражения со вторыми производными выполняется оценка:
\begin{equation}\label{ots-1-2}
\begin{array}{c}
\Biggl| \displaystyle\frac{1}{2}\displaystyle
\int\limits_{0}^{t}\Bigl[\dfrac{1}{(\varepsilon\sqrt{2\pi})^{n}}\displaystyle\int\limits_{\mathbb{R}^{n} }\prod\limits_{i=1}^{n}\exp\left\{-\frac{(y_{i}-x_{i}(\tau))^{2}}{2\varepsilon^{2}}   \right\}
 \frac{\partial^{\,2}}{\partial y_{i}\partial y_{j}}F(\tau,{\bf y}) d\Gamma({\bf y}) -\\
-\displaystyle\frac{\partial ^{2} F(\tau;{\bf x})}{\partial x_{i} \partial x_{j} }\bigl|_{{\bf x}={\bf x}(\tau) }\Bigr]{b}_{i\,k}(\tau){b}_{j\,k}(\tau)d\tau \Biggr|\leq
\\
\leq
\dfrac{\varepsilon^{n}2^{n+1}}{(\sqrt{2\pi})^{n}}
\displaystyle\int\limits_{0}^{t}\Bigl|{b}_{i\,k}(\tau){b}_{j\,k}(\tau)\Bigr|
L^{(8)}(\tau)d\tau.
\end{array}
\end{equation}

Поскольку $L^{(s)}(t)$, $s=1,3-6,8$ -- неслучайные ограниченные функции, функции ${a}_{i}(t)$ и ${b}_{j\,k}(t)$ удовлетворяют условиям \eqref{GrindEQ__2_1_2_}, то при $\varepsilon\downarrow 0$ выражения в \eqref{ots-Q}, \eqref{ots-der-a}, \eqref{ots-Fg}--\eqref{ots-der-D}, \eqref{ots-1-2} не превосходят 0.

Интеграл по винеровскому процессу $\displaystyle\int\limits_{0}^{t}f(\tau)dw(\tau)$ определен для функций $f(\tau)\in H_{2}[0;t]$, таких, что с вероятностью единица выполняется условие $\displaystyle\int\limits_{0}^{t}f^{2}(\tau)d \tau < \infty$ (\cite[с.~12]{GS_68}), т. е. выполняются ограничения \eqref{GrindEQ__2_1_2_}, и следовательно,  выражения в \eqref{ots-D} и  \eqref{ots-der-F} при $\varepsilon\downarrow 0$ также не превосходят 0.

Следовательно, при совпадении начальных условий имеем:
\begin{equation}\label{itog-ots}
\lim\limits_{\varepsilon\downarrow 0}\Bigl|F_{\varepsilon}(t;{\bf x}(t))-F(t;{\bf x}(t))\Bigr|\leq 0.
\end{equation}
Неравенство \eqref{itog-ots} имеет смысл только в виде
\begin{equation*}\label{itog-ots}
\lim\limits_{\varepsilon\downarrow 0}\Bigl|F_{\varepsilon}(t;{\bf x}(t))-F(t;{\bf x}(t))\Bigr|= 0.
\end{equation*}
Это означает выполнение равенств \eqref{M-ots} и \eqref{r01}:
\begin{equation*}\label{d01}
F(t;{\bf x}(t))= {\underset{\varepsilon \downarrow 0}{\mbox { \textrm{l.i.m.}}}} \,
F_{\varepsilon}(t;{\bf x}(t)).
\end{equation*}
Следовательно, имеет место утверждение
\eqref{GrindEQ__2_5_2_} и теорема доказана.
\end{proofs}

Отметим, что обобщенная формула Ито -- Вентцеля включает в себя, как частные случаи, рассмотренные выше правила построения дифференциалов \eqref{odu-prav} и \eqref{GrindEQ__1_4_2_}. Таким образом, ее можно использовать и для детерминированных уравнений.

\section*{Заключение}

Обобщенная формула Ито -- Вентцеля для НЦПМ требует более жестких ограничений на класс весовых функций при пуассоновских возмущениях по сравнению с аналогичной формулой для случая ЦПМ (см., например, \cite{D_02}). Применение классических методов теории стохастических дифференциальных уравнений на основе аппроксимации специального вида  $\delta$-функции позволили определить  достаточные условия существования и единственности решения задачи Коши для ОСДУ Ито с  НЦПМ,  представляющего собой обобщенную формулу Ито -- Вентцеля для нецентрированной меры.

Отметим, что для построения интегралов по ЦПМ достаточно, чтобы с вероятностью единица выполнялось условие:
$
\displaystyle\int |f(\gamma )| ^{2} \Pi (d\gamma )<\infty,
$
где
$\Pi(\mathcal{A})$ -- мера множества $\mathcal{A}$ из $\sigma$-алгебры борелевских множеств в $\mathbb{R}^{n'}$.
При работе с НЦПМ добавляется еще условие:
$
\displaystyle\int |f(\gamma )| \Pi (d\gamma )<\infty
$
(см. \cite[c.\,255;\,c.\,254]{GS_68}).

 \textit{Выражаю благодарность проф. В. А. Дубко за предложенную для доказательства идею.}

\end{document}